\documentclass[12pt]{article}
\usepackage{amssymb}
\usepackage{color}
\usepackage{geometry}

\newtheorem{Assumption}{Assumption}
\newtheorem{Theorem}{Theorem}
\newtheorem{Proposition}{Proposition}
\newtheorem{Remark}{Remark}

\usepackage[english]{babel}
\usepackage[latin1]{inputenc}
\usepackage{times}
\usepackage[T1]{fontenc}
\title
{ {\normalsize\tt\hfill\jobname.tex}\\
On Poisson equations with a potential in the whole space for "ergodic" generators
 \\~\\
Alexander Veretennikov% <-this % stops a space
\footnote{ University of Leeds, UK; National Research University Higher School of Economics, and Institute for Information Transmission Problems, Moscow, Russia, email: a.veretennikov @ leeds.ac.uk
}
}

\begin{document}

\maketitle

\begin{abstract}
{\small  
In \cite{PV01, PV03, PV05, Ver11, KV11}
Poisson equation \emph{in the whole space} 
was studied for so called ergodic generators $L$ 
corresponding to homogeneous Markov diffusions 
($X_t, \, t\ge 0$) in $\mathbb R^d$.  Solving this equation is one of the main tools for  
{\it diffusion approximation} in the theory of stochastic averaging and homogenisation. Here a similar equation {\it with a potential} is considered, 
firstly because it is natural for PDEs, 
and secondly with a hope that it may be also 
useful for some extensions related to homogenization and averaging. 
%Let $\mu$ be the (unique) invariant probability measure of  $X$. 
}
\end{abstract}

\section{Introduction}
Let us consider a stochastic differential equation in $\mathbb R^d$, 
\begin{equation}\label{eq1}
dX_t = \sigma(X_t)\,dB_t + b(X_t)\,dt, \quad X_0=x. 
\end{equation} 
Assume the matrix function $\sigma$ and the vector function $b$ are  Borel bounded, 
%and the vector function $b$ is Borel with a linear growth condition, 
and $\displaystyle a = (a^{ij}(x), \, x\in \mathbb R^d):= \frac12 \sigma\sigma^*$ uniformly non-degenerate; some further conditions will be assumed in the sequel (here $\sigma^*$ stands for the transposed matrix $\sigma$), however, the above already suffices for the existence of solution \cite{Kry69}.
% {\color{red}(to check linear growth!..)}. 
Denote 
\[
L = \sum_{i,j=1}^{d} a^{ij}(x) \frac{\partial ^2}{\partial x^i \partial x^j} + \sum_{i=1}^{d} b^{i}(x) \frac{\partial}{\partial x^i}. 
\]
The Poisson equation 
\begin{equation}\label{eq01}
Lu(x) = -f(x), \quad x\in R^d,
\end{equation}
is one of the well-known equations in mathematical physics. Only relatively recently, in the last two decades it was understood how useful is this equation  {\em without boundary conditions and in the whole space $\mathbb R^d$}: namely, it is a powerful tool in the theory of averaging and homogenization, see \cite{Bahvalov, Ethier, PSV}, et al. This understanding was the reason for the series of papers \cite{PV01, PV03, PV05} and for some further versions and extensions in \cite{KV11,  Ver12}. In all papers in the latter references such equations were {\em without zero order terms} (also known as potentials). On the other hand, equations {\em with} potentials are also very frequent in physics and even more popular than without them. In the cases where the author presented these results at the PDE seminars, the most frequent question was, why zero order terms are not included in the equation? 

Hence, the goal of this paper is to transfer some of the recent advances about Poisson equation ``in the whole space'' without a potential to the equation 
\begin{equation}\label{eq02}
Lu -cu= -f, 
\end{equation}
with a potential $c$. Note straight away that the case $\inf_x c(x) > 0$ is the most simple one where -- at least, for bounded functions $f$ -- convergence of the integral representing the solution
\begin{equation}\label{repsol}
u(x) = \int_0^\infty {\mathbb E}_x \exp(- \int_0^t c(X_s)\,ds)f(X_t)\,dt,
\end{equation}
and the equation itself for this representation follow automatically. So, we will concentrate on the more interesting situations where either $c$ is {\em not} separated away from zero yet remaining non-negative, or even if the function $c$ may change its sign, with a hope that in the future it could be possibly useful, in particular, for controlled Markov processes and, perhaps, for a probabilistic interpretation of the Helmholtz equation.  The problem of equations with parameters is not addressed here.  

The paper consists of Introduction; Reminder of equations without potentials; Main results; Proofs; and the latter part is additionally split into several sections and subsections. 

\section{Equation without potential: quick reminder} %(\ref{eq01})}
Let us present the main assumptions. They will be used in the next section, too, except that the first one (A1) will be replaced by a stronger assumption (A5). Also note that the centering condition (A4) is needed only in this section, and also in one case out of three cases in the Main results. 

\begin{itemize}

\item[(\bf A1)]
\begin{Assumption}[polynomial  recurrence]\label{As1}

\begin{equation}\label{poly_ass}
\limsup_{|x|\to\infty}\langle b(x), x\rangle = - \infty.
\end{equation}

\end{Assumption}

%\item[(\bf A3)]
%\begin{Assumption}[A3 -- non-degeneracy]\label{As3}
%$a:=\sigma\sigma^*$ is uniformly non-degenerate, bounded 
%\end{Assumption}

\item[(\bf A2)]
\begin{Assumption}[boundedness and non-degeneracy]\label{As2}
~\\
The coefficient $\sigma\sigma^*$ is uniformly non-degenerate
\begin{equation}\label{snd}
\inf_{x} \inf_{|\xi|=1} \xi^*\sigma\sigma^*(x)\xi >0;
\end{equation}
$\sigma$, $b$, $f$ and $c$  are Borel bounded.
%; drift $b$ satisfies a linear growth condition
%\begin{equation}\label{blg}
%|b(x)| \le C(1+|x|).
%\end{equation}
\end{Assumption}

\item[(\bf A3)]
\begin{Assumption}[regularity]\label{As3}
~\\
The coefficient $a=\sigma\sigma^*$ is uniformly continuous.
% {\small (plus something else, see later)}.  
\end{Assumption}

\item[(\bf A4)]
\begin{Assumption}[$f$--centering]\label{As4}
\begin{equation}\label{as-centre}
 \int f(x)\mu(dx) = 0.
\end{equation}

\end{Assumption}

\item[(\bf A5)]
\begin{Assumption}[ exponential recurrence]\label{As5}
\hspace{3.4cm}
\begin{equation}\label{ld_ass}
\limsup_{|x|\to\infty}\langle b(x), x/|x|\rangle = - r < 0. 
\end{equation}
\end{Assumption}

\end{itemize}
\begin{Remark}
The assumption (A3) may be totally  dropped in one-dimensional case. In the multi-variate case it is assumed so as to guarantee weak uniqueness of solutions of the equation (\ref{eq1}). The assumption (A5) may be strengthened to $\limsup_{x}\langle b(x), x/|x|\rangle = - \infty$, in which case some other references on large deviation results would be required but the results would be established not just for small values of $\epsilon$ only; however, this would also require new references about convergence rates and mixing because for this assumption to be non-empty, the drift may not be globally bounded. The boundedness of other coefficients may also be relaxed, but we do not pursue this goal here,
\end{Remark}

Here is the main result from \cite{PV01} about the equation (\ref{eq01}); we do not show further advances already established, which relate to the more involved case of equations with parameters.  Note that under the assumptions (A1) -- (A3), the equation (\ref{eq1}) has a weakly unique Markov and strong Markov solution $(X_t)$ with a unique stationary measure $\mu$ (cf. \cite{Kry69, Ver00}).

\begin{Proposition}[\cite{PV01}]
Under the assumptions (A1) -- (A4),  the equation (\ref{eq01}) has a solution $u$ in  $\bigcap_{p>1} W^{2}_{p, loc}$. This solution is itself centered, it has a moderate growth (i.e., not faster than some polynomial), and this solution is unique in this class of functions. The representation (\ref{repsol}) with $c\equiv 0$ holds true for this solution. 
\end{Proposition}

Recall that the assumption (A5) is not needed in this Proposition. On the other hand, where it will be used (in the next sections), it clearly replaces the weaker assumption (A1).

\section{Equation with a potential}
Now we turn to the main goal of the paper, the equation (\ref{eq02}) with a potential $c$. As it was already mentioned, a natural candidate for the solution is the function
\[\displaystyle 
 u(x) = \int_0^\infty \mathbb E_x \exp\left(-\int_0^t c(X_s)\,ds
\right) f(X_t)\,dt, 
\] 
see (\ref{repsol}), 
of course, provided that this expression is well-defined. Recall that in the sequel it will be assumed that both $f$ 
and $c$ are bounded. Beside the most simple case $\inf_x c(x) > 0$, we are able to tackle three different situations. In all of these cases we will assume the assumption (A5), which, actually, replaces the weaker one (A1),  see below. 

\bigskip

{\bf Case 1}: $c(x) = \epsilon  c_1(x)$, where $\epsilon > 0$ is \emph{small enough}, and $\int c_1 d\mu > 0$.

\medskip

{\bf Case 2}: $c(x) \ge 0, \; \& \; \int c\, d\mu > 0$. {\it It is not assumed that (\(c\) is small or bounded away from zero here).}

\medskip

{\bf Case 3}: $c(x) = -\epsilon < 0$ with $\epsilon$ -- which is a constant -- small enough.

\bigskip

\noindent
Note that all three cases do not include each other, although the case 1 and the case 2 do intersect. 
In all three cases we assume (A2) -- (A5).
The question about the case 2 was suggested by A. Piunosvkiy; hopefully, it might be useful in the theory of  controlled Markov processes (cf. with \cite{Piu}). 

~

\noindent
%In the case 2 some further considerable extension is possible. 
In the cases 1 and 2, denote 
$$
\bar c = \int c d\mu, 
$$
and in the case 1 also 
$$
\bar c_1 = \int c_1 d\mu.
$$
Clearly, $\bar c = \epsilon \bar c_1$. 

\begin{Theorem}\label{thm1}
In the cases 1 -- 2, under the assuptions (A2)--(A3) and (A5), the function $(u(x), \, x\in \mathbb R^d)$ given by the representation (\ref{repsol}) is a continuous solution in the Sobolev classes $W^{2}_{p, loc}$ for each $p>0$ of the Poisson equation (\ref{eq02}). This solution admits the bounds 
\begin{equation}\label{ubds}
|u(x)| \le C \exp(\gamma |x|), 
\end{equation}
with any $\gamma>0$ and a corresponding $C = C(\gamma)$. 
In the case 3 the same assertions hold true under the assumptions (A2) -- (A5), and, in addition, the function given by the formula (\ref{repsol}) is centered.
\end{Theorem}
Note that in the cases 1--2 there is no need for the centering assumption (A4). The proof is split into the next three section devoted to convergence (the section 4), existence of derivatives and verification of the equation (the section 5); uniqueness (the section 6).

\section{Proof: convergence}
In all cases we will use the bound from \cite{Ver87}
\begin{equation}\label{expmix}
\|Q_t(x,dy) - \mu(dy)\|_{TV} \le C \exp(\gamma |x|) \exp(-\lambda t), \quad t\ge 0, 
\end{equation}
where $Q_t(x,dy)$ is the transition kernel of the process $X_t$, and $\mu$ is its unique invariant measure, and ``TV'' is the total variation distance for two measures. Note that the inequality (\ref{expmix}) may be read as follows: there exists $\lambda_0>0$ such that for any $\lambda \in (0, \lambda_0)$, there exists $\gamma>0$ such that (\ref{expmix}) holds; yet, it may be also read as follows: there exists $\gamma_0>0$ such that for any $\gamma \in (0, \gamma_0)$, there exists $\lambda>0$ such that (\ref{expmix}) holds true. Several close but a little different corollaries from this inequality will be used in the sequel.

~

\noindent
Another bound from \cite{Ver87} reads: there exists $\gamma_0>0$ such that  for any $\gamma \in (0,\gamma_0)$ there exists $C>0$ such that  
\begin{equation}\label{expmom}
\sup_{t\ge 0} \, \mathbb E _{x}\exp(\gamma |X_t|) \le C \exp(\gamma |x|). 
\end{equation}

\subsection{Case 1, locally uniform convergence }
In this section we show convergence of the integral, 
\[
 u(x) := \int_0^\infty \mathbb E_x \exp\left(-\int_0^t c(X_s)\,ds
\right) f(X_t)\,dt,
\]
in the case with 
$
c(x) = \epsilon  c_1(x)$, 
$\bar c_1=\int c_1 d\mu > 0,$ with $\epsilon >0 $ small enough.  
Denote 
$$
H_T(\beta,x) := T^{-1} \ln \mathbb E_x \exp 
\left(\beta \int_0^T c_1(X_s)\,ds\right), 
$$ 
or, equivalently, 
\[
\mathbb E_x \exp \left(\beta \int_0^T c_1(X_s)\,ds\right) = \mathbb E_x \exp(T\,H_T(\beta,x)). 
\]
Let 
$$
H(\beta) := \lim_{T\to\infty} H_T(\beta,x), \qquad \beta\in 
R^d.
$$
Under the condition (A5), this limit does exist for all values of \(\beta\) locally uniformly with respect to $\beta$ and $x$, see, e.g., \cite[Theorem 1]{Ver95}. Since $H_T(\beta,x)$ is convex with respect to $\beta$, and due to this locally uniform convergence, the limiting function $H(\beta)$ is also convex. Since $H$ is clearly finite for any value of $\beta$, it is also continuous as any finite convex function. It follows that $H$ is differentiable at the origin $\beta = 0$, see, e.g., \cite{Ver-viln}. Because the reference may be not very well accessible, we recall the idea of this simple reasoning. 
For this aim it is convenient to perform the following transformation,  
$$
H^{1}_T(\beta,x):= T^{-1} \ln \mathbb E_x \exp 
\left(\beta \int_0^T (c_1(X_s) - \bar c_1)\,ds\right) 
\equiv H_T(\beta,x) - \beta c_1. 
$$
Further, due to the Law of Large Numbers -- by virtue of a good mixing for \(X\), see, e.g., \cite{Ver00} -- we have,  
$$
({\mathbb P}_x)\lim_{T\to\infty} T^{-1}  \int_0^T (c_1(X_s) - \bar c_1)\,ds = 
0, 
$$  
and also $
\lim_{T\to\infty} T^{-1} {\mathbb E}_x \int_0^T c_1(X_s)\,ds = 
\bar c_1$. Moreover, under the assumption (A5) due to the exponential mixing bound (cf. \cite{Ver87, Ver95}), for any $\epsilon>0$ there exist $C=C(x), \lambda>0$, such that an exponential 
%Bienaim\'e--Chebyshev--Markov 
inequality holds, 
\begin{equation}\label{bcm1}
{\mathbb P}_x (| T^{-1}  \int_0^T (c_1(X_s) - \bar c_1)\,ds|\ge \epsilon) 
\le  C(x)\exp(-\lambda t),  
\end{equation}
with some $C(x) = \exp(\gamma |x|), \; \gamma>0$.
Note that here $\lambda>0$ does not depend on $x$. With $\epsilon>0$ small enough, and using the split of unity 
$$
1 = 1(|T^{-1}  \int_0^T (c_1(X_s) - \bar c_1)\,ds| <\epsilon) 
+ 1(|T^{-1}  \int_0^T (c_1(X_s) - \bar c_1)\,ds| \ge \epsilon), 
$$
and the elementary inequality $a+b \le 2 \, (a\vee b)$ and, hence, $\ln (a+b) \le \ln (2 \, (a\vee b))$ (with $a,b>0$), 
we now compute for any $|\beta|\le b$, say, 
\begin{eqnarray*}
&\displaystyle H^{1}_T(\beta,x) = T^{-1} \ln \left(\mathbb E_x \exp 
(\beta \int_0^T (c_1(X_s) - \bar c_1)\,ds) 1(|T^{-1}  \int_0^T (c_1(X_s) - \bar c_1)\,ds| <\epsilon) \right.
 \\\\ &\displaystyle \left. 
 + \mathbb E_x \exp 
\left(\beta \int_0^T (c_1(X_s) - \bar c_1)\,ds\right) 1(|T^{-1}  \int_0^T (c_1(X_s) - \bar c_1)\,ds| \ge \epsilon) \right)
 \\\\
&\displaystyle \le T^{-1} \ln \left(2 (\exp 
(|\beta| \,T \epsilon) \vee \left(\exp 
(|\beta| \,T \|c_1-\bar c_1\, \|\right)C(x)\exp(-\lambda T))\right)
 %\\\\
\le 2 b\epsilon.
\end{eqnarray*}
Here $b>0$ may be taken small enough in comparison to $\lambda$ (since the latter does not depend on $\beta$). This implies that for a fixed $x$ -- and, actually, locally uniformly with respect to $x$ -- the function $H^1_T(\beta, x) = o(|\beta|)$ uniformly in $T\to\infty$. Thus, $H^1(\beta):= H(\beta) - \bar c_1\beta = o(|\beta|)$, which, clearly, means that $H(\beta)$ is differentiable at zero and that $H'(0) = \bar c_1$. Note that this is also in accordance with the fact that $H_
T(\beta,x) \to H(\beta)$ and  since both functions are 
convex in \(\beta\), we also have $H'_T(0,x) \to H'(0)$
(see, e.g., \cite{Rock}). For any $\delta>0$ we may assume that 
$$
|H_T(\beta,x)  - H(\beta)| \le \delta + o_T(1), \quad \mbox{as} \quad T\to\infty.
$$ 
In any case, since $H'(0) = \bar c_1>0$, in some neighbourhood of zero, 
\begin{equation}\label{Hmonotone}
H(\beta)>0, \quad \beta>0, \qquad \& \qquad H(\beta)<0, 
\quad \beta<0.
\end{equation}
Therefore, convergence of the integral in the definition of $u$ follows from 
(\ref{Hmonotone}). Indeed, we estimate, for $\epsilon>0$ small and independent on $x$,  and with any $\delta>0$ and taking $\beta = -\epsilon$ we have, 
\begin{eqnarray*}
& \displaystyle |u(x)| \le \|f\|_B \, \int_0^\infty \mathbb E_x \exp\left(-\varepsilon \int_0^t c_1(X_s)\,ds\right) \,dt
 %\\\\
= \|f\|_B \, \int_0^\infty \exp(t H_t(-\epsilon,x) \,dt 
 \\\\
& \displaystyle \le  \|f\|_B \, \int_0^\infty \exp(t H(-\epsilon) +  \delta +  o_t(1) + \gamma |x|)\,dt 
< \infty. 
\end{eqnarray*}
Also note that here $\gamma>0$ may be chosen arbitrarily small, which means that the rate of growth of the function $u$ is slower than any exponential of $|x|$. (In fact,  some better polynomial growth bound on $u(x)$ holds true, too.) More precisely, for any $\gamma>0$ small enough (and, hence, actually, for any $\gamma>0$) there exists $C>0$ such that 
\begin{equation}\label{ubd1}
|u(x)| \le C \exp(\gamma |x|).
\end{equation}

\subsection{Case 2, locally uniform convergence }

Recall that we wish to show convergence of the integral,
 
\[
 u(x) := \int_0^\infty\mathbb  E_x \exp\left(-\int_0^t c(X_s)\,ds
\right) f(X_t)\,dt,
\]
in the case where  
$$
c(x) \ge 0, \quad \& \quad \int c\,d\mu >0.
$$
We will use an exponential estimate similar to (\ref{bcm1}),
\begin{equation}\label{e-expbd}
P_x(\int_0^t c(X_s)\,ds < (\bar c - \delta)t) \le C\exp(
\nu|x|)\exp(-\lambda t),
\end{equation}
along with a split of unity into two indicators. It follows 
that again the expression for $u$ is well defined.
Indeed, 
\begin{eqnarray*}
&\displaystyle \int_0^\infty \mathbb E_x \exp\left(-\int_0^t c(X_s)\,ds \right) f(X_t)\,dt 
 \\\\
&\displaystyle = \int_0^\infty \mathbb E_x \exp\left(-\int_0^t c(X_s)\, \,ds \right) f(X_t)
%\\\\ \times 
\, 1(\int_0^t c(X_s)\,ds < (\bar c - \delta)t)\,dt
 \\\\
&\displaystyle +\int_0^\infty \mathbb E_x \exp\left(-\int_0^t c(X_s)\, \,ds \right) f(X_t)
%\\\\ \times 
1(\int_0^t c(X_s)\,ds \ge  (\bar c - \delta)t)\,dt. 
\end{eqnarray*}
Here the second term clearly converges for small \(\delta\), while the first term converges due to the assumption \(c\ge 0\) and \(f\) bounded, because of the inequality (\ref{e-expbd}), as required. 

~

Also, with the help of (\ref{expmom}) it follows, 
\begin{equation}\label{ubd2}
|u(x)| \le C \exp(\gamma |x|).
\end{equation}

\subsection{Case 3, locally uniform convergence and centering }
Convergence along with the bound on $|u(x)|$ in this case follows straightforward from the inequality (\ref{expmix}), 
if we only admit that $-c \equiv \epsilon < \lambda$. The bound on $|u(x)|$ then reads, 
\begin{equation}\label{ubd3}
|u(x)| = |\int_0^\infty \mathbb E_x \exp\left(\epsilon \, t\right) f(X_t)\,dt | \le 
C \|f\|_B(\lambda-\epsilon)^{-1}\exp(\gamma |x|).  
\end{equation}

~

The centering condition holds true due to the same centering assumption on $f$: by virtue of Fubini's theorem, 

\begin{eqnarray*}
&\displaystyle \int u(x) \mu(dx) = \int \mu(dx) \int_0^\infty\mathbb  E_x \exp\left(+\epsilon \, t\right) f(X_t)\,dt 
 \\\\
&\displaystyle =  \int_0^\infty \int \mu(dx) \mathbb E_x \exp\left(+\epsilon \, t\right) f(X_t)\,dt =
0,
\end{eqnarray*}
where the fact was used that the measure $\mu$ is stationary and, hence, $\int \mu(dx) \mathbb E_x f(X_t) = 0$ for each $t\ge 0$. 

\section{Proof: other properties}

\subsection{Verification of the equation: simplified version}

As we already know from the previous section, the function $u$ given by the representation (\ref{repsol}) is well-defined, that is, the integral in the right hand side converges for all $x \in \mathbb R^d$. Let us argue why $u$ is, indeed, a solution of the Poisson equation. To make explicit the idea, assume for simplicity continuity of both functions  
$c$ and $f$ and suppose that existence of two (classical) derivatives of the function $u$ is known; later on it will be shown how to drop all these additional assumptions, including classical derivatves instead of Sobolev ones. By the Markov property, 
\begin{eqnarray*}
&\displaystyle  u(x) = \int_0^T \mathbb E_x \exp\left(-\int_0^t c(X_s)\,ds\right) 
f(X_t)\,dt 
 \\
&\displaystyle + 
\mathbb  E_x \exp\left(-\int_0^T c(X_s)\,ds\right)\, u(X_T),
\end{eqnarray*}
from which,
\begin{eqnarray*}
&\displaystyle  -f(x) = -\lim_{T\to 0} T^{-1}\, \int_0^T \mathbb E_x \exp\left(-\int_0^t c(X_s)\,ds\right) 
f(X_t)\,dt 
 \\\\
&\displaystyle = \lim_{T\to 0} T^{-1} \left(\mathbb E_x u(X_T)\exp\left(-
\int_0^T c(X_s)\,ds\right) - u(x)\right)
 \\ \\
&\displaystyle = \lim_{T\to 0} \frac1{T} \mathbb E_x \int_0^T e^{\displaystyle -\int_0^t c(X_s)ds}
(Lu - cu)(X_t)\,dt = Lu(x) - cu(x), %\hspace{4cm}
\end{eqnarray*}
as required. However, as we said earlier, in the sequel we aim to justify the equation without the additional assumption about continuity.

\subsection{Continuity of solution $u$}
This continuity of $u$ will be used in the proof of existence of two Sobolev derivatives in the next subsection. 
Actually, we shall see a bit more than just continuity:   in all three cases 1 -- 3 it will be shown  that the integral for $u$ converges to a continuous limit  \emph{locally uniformly} with respect to $x$. 
So, similarly to \cite{PV01} we obtain, 
\begin{eqnarray*}
& \displaystyle u(x) = \lim_{N\to\infty} u^N(x) %\hspace{2cm}
 \\ \\
& \displaystyle:= \lim_{N\to\infty}\mathbb  E_x \int_0^N  \exp\left(-\int_0^t c(X_s)
\,ds 
\right) f(X_t)\,dt,
\end{eqnarray*}
where 
%by Fubini's theorem \(\displaystyle E\int_0^N \ldots = \int_0^N E \ldots \), and 
$u^N(x)$ is continuous as a solution of the Cauchy 
problem for a parabolic differential equation, see \cite{KS}. 
So, the limit is also continuous, due to the locally uniform convergence. 
{\it Note that neither continuity of $f$ nor of $c$ was used in this consideration.}

\subsection{\bf Two Sobolev derivatives for $u$ and verification of the equation
%: cf. [P, V, 2001]
}

Consider $\tau:=\inf(t:\, X_t\not\in B)$ and the following equation in 
the ball $B=B_1(x_0)\equiv \{x\in \mathbb R ^d:\; |x-x_0|\le 1\}$,
$$
Lv - cv = -f \quad \mbox{in} \quad B, \qquad v|_{\partial B}
 = u
$$
Since we already know that $u$ is continuous, this boundary condition is well-defined. There is a unique   
solution $v \in \bigcap_{p\ge 1} W^{2}_{p}$ in $B$ 
%[Ladyzhenskaya, Ural'tseva] 
\cite{LSU}, which by virtue of It\^o--Krylov's formula \cite{Kry77} admits a representation, 
\begin{eqnarray*}
 v(x) = \mathbb E_x\left(\int_0^{\tau} \exp\left(-\int_0^t c(X_s) \,
ds\right) f(X_t)\,dt 
 %\right. \\ \left. 
+ \exp\left(-\int_0^{\tau} c(X_s)\,ds\right) u(X_{\tau})
\right). 
\end{eqnarray*}
Due to the strong Markov property, exactly the same representation 
holds true for $u(x)$ in the left hand side; so,  $u \equiv v$ on $B$, i.e.,  \begin{eqnarray*}
 u(x) =\mathbb E_x\left(\int_0^{\tau} \exp\left(-\int_0^t c(X_s) \,
ds\right) f(X_t)\,dt 
 %\right. \\ \left. 
+ \exp\left(-\int_0^{\tau} c(X_s)\,ds\right) u(X_{\tau})
\right). 
\end{eqnarray*}
%Indeed, ...
Hence, $u \in \bigcap_{p\ge 1} W^{2}_{p}(B) \subset C(B)$ (see \cite{LSU}). This consideration also justifies the equation for $u$ without the additional assumption about continuity of $f$ and $c$.

\section{Proof: uniqueness of solution}

\subsection{Uniqueness, case 1}

Uniqueness may be shown in a standard manner for the class 
of functions satisfying the moderate growth established 
earlier. In all three cases the calculus is the same. For the difference of two solutions $v = u^1 - 
u^2$ we have $Lv - cv = 0$. So, using moment inequalities 
and a standard localization procedure,  by applying It\^o--Krylov's 
formula and taking expectations, we get 
$$
v(x) = u^1(x) - u^2(x) =\mathbb  E_x \exp(-\int_0^t c(X_s)\,ds) v(X_
t).
$$
We now use a unity split 
$$1 = 1(\int_0^t c_1(X_s)\,ds 
\ge (\bar c_1 - \delta)t) + 1(\int_0^t c_1(X_s)\,ds < (\bar 
c_1 - \delta)t)$$ and an exponential estimate 
$$
\mathbb P_x(\int_0^t c_1(X_s)\,ds < (\bar c_1 - \delta)t) \le C\exp(
\nu|x|)\exp(-\lambda t).
$$
Then it follows that $v(x) \equiv 0$. Indeed, due to the bound (\ref{bcm1}), we estimate with any $t>0$,  
\begin{eqnarray*}
&\displaystyle  |v(x)| \le \mathbb E_x \exp(-\int_0^t c(X_s)\,ds) |v(X_t)| 
 \\ \\
&\displaystyle = \mathbb E_x \exp(-\int_0^t c(X_s)\,ds) |v(X_t)|1(\int_0^t c(X_s)\,
ds \ge (\bar c - \delta)t)
 \\\\
&\displaystyle  + \mathbb E_x \exp(-\int_0^t c(X_s)\,ds) |v(X_t)|1(\int_0^t c_1(X_
s)\,ds < (\bar c_1 - \delta)t)
 \\\\
&\displaystyle \le  \mathbb E_x \exp(-\epsilon (\bar c_1 - \delta)t) C \exp(
\nu|X_t|) 
 \\\\
&\displaystyle + \exp(\epsilon \|c_1\| t) 
 %\\\\ \times 
\, (\mathbb E_x|v(X_t)|^2)^{1/2} (\mathbb P_x(\int_0^t c_1(X_s)\,ds < (
\bar c_1 - \delta)t))^{1/2}
 \\\\
&\displaystyle \le C \exp(-\epsilon (\bar c_1 - \delta)t) C \exp(\nu 
|x|) 
 \\\\
&\displaystyle  + C \exp(\nu' |x|) \exp(-(\lambda/2 - \epsilon \|c_1\|)
t) \to 0, \quad t\to\infty.
\end{eqnarray*}
In teh middle of the calculus we have applied Caushy--Bouniakovsky--Schwarz' inequality. So, if $\epsilon>0$ was chosen small enough, it shows that $u^1 \equiv u^2$, which completes the proof of the Theorem \ref{thm1} in the case 1. 

~

%There is a natural guess that the cases 1 and 2 might be combined in the way where \(c\) may change sign, however, its only {\it positive part} is small:
%\[
%c(x) = \epsilon c^1_-(x) + c^1_+(x)_, \quad \& \quad c^1_+ \not \equiv 0. 
%\] 
%We do not pursue this goal in this paper, although, this seems reasonable and not too involved. Apparently, a completely open question is whether or not the case $c^1_+ \equiv 0$ with $c^1_-$ {\em variable}  may be tackled. 

\subsection{Uniqueness, case 2}

For the difference of two solutions $v = u^1 - 
u^2$ we have $Lv - cv = 0$. So, using moment inequalities 
and a standard localization procedure,  by applying It\^o--Krylov's 
formula and taking expectations, we get 
$$
v(x) = u^1(x) - u^2(x) = \mathbb E_x \exp(-\int_0^t c(X_s)\,ds) v(X_
t).
$$
We now use a unity split 
$$1 = 1(\int_0^t c(X_s)\,ds 
\ge (\bar c - \delta)t) + 1(\int_0^t  c(X_s)\,ds < (\bar 
c_1 - \delta)t)$$ and an exponential estimate 
$$
\mathbb P_x(\int_0^t c(X_s)\,ds < (\bar c - \delta)t) \le C\exp(
\nu|x|)\exp(-\lambda t).
$$
Then it follows that $v(x) \equiv 0$. Indeed, recall that in the case 2, $c\ge 0$, Hence,  due to the bound (\ref{e-expbd}), we estimate, with any $t>0$,  
\begin{eqnarray*}
&\displaystyle  |v(x)| \le \mathbb E_x \exp(-\int_0^t c(X_s)\,ds) |v(X_t)| 
 \\ \\
&\displaystyle = \mathbb E_x \exp(-\int_0^t c(X_s)\,ds) |v(X_t)|1(\int_0^t c(X_s)\,
ds \ge (\bar c - \delta)t)
 \\\\
&\displaystyle  + \mathbb E_x \exp(-\int_0^t c(X_s)\,ds) |v(X_t)|1(\int_0^t c(X_
s)\,ds < (\bar c - \delta)t)
 \\\\
&\displaystyle \le  \mathbb E_x \exp(- (\bar c - \delta)t) C \exp(
\nu|X_t|) 
 \\\\
&\displaystyle + (\mathbb E_x|v(X_t)|^2)^{1/2} (\mathbb P_x(\int_0^t c(X_s)\,ds < (
\bar c - \delta)t))^{1/2}
 \\\\
&\displaystyle \le C \exp(-  (\bar c - \delta)t) C \exp(\nu 
|x|) 
 \\\\
&\displaystyle  + C \exp(\nu' |x|) \exp(- \lambda/2)
t) \to 0, \quad t\to\infty.
\end{eqnarray*}
So, $u^1 \equiv u^2$, which completes the proof of the Theorem \ref{thm1} in the case 2. 

~

\subsection{Uniqueness, case 3}

For the difference of two (centered) solutions $v = u^1 - 
u^2$ we have $Lv - cv = 0$. So, using moment inequalities 
and a standard localization procedure,  by applying It\^o--Krylov's 
formula and taking expectations, we get 
$$
v(x) = u^1(x) - u^2(x) = \mathbb E_x \exp(-\int_0^t c(X_s)\,ds) v(X_
t) \equiv \exp(+ \epsilon t) \mathbb E_x  v(X_t).
$$
Recall the bound $|v(x)| \le C\exp(\gamma |x|)$, where $\gamma>0$ can be chosen arbitrarily small (and, of course, respectively, $C$ depends on $\gamma$). By using the bounds (\ref{expmix}) and (\ref{expmom}) and taking $0< \delta<\bar c$, and due to the centering property of $v$, we estimate,

\begin{eqnarray*}
&\displaystyle 
 |v(x)| \le \exp(+ \epsilon t) |\mathbb E_x v(X_t)| = \exp(+ \epsilon t) \,
|\int v(y) Q_t(x,dy)|
 \\ \\
&\displaystyle = \exp(+ \epsilon t) \, \left|\int v(y) (Q_t(x,dy)-\mu(dy)\right|
 \\\\
&\displaystyle  \le \exp(+ \epsilon t) \, \left(\int v^2(y) (Q_t(x,dy)+\mu(dy)\right)^{1/2} \left(\int |Q_t(x,dy)-\mu(dy)|\right)^{1/2} 
 \\\\
&\displaystyle  \le \exp(+ \epsilon t) \, C\exp(\gamma |x|)\exp(-t\,\lambda/2) \to 0, \quad t\to\infty, 
\end{eqnarray*}
if $\epsilon>0$ is small enough. We used, in particular, Cauchy--Bouniakovsky--Schwarz' inequality and the fact that $\mu$ integrates exponentials $\exp(\gamma|x|)$ with small $\gamma$. 
So, $u^1 \equiv u^2$, which completes the proof of the Theorem \ref{thm1} in the case 3. 

~

\section*{Acknowledgements}
This work was prepared within the framework of a subsidy granted to the HSE by the Government of the Russian Federation for the implementation of the Global Competitiveness Program, and supported by the RFBR grant 14-01-00319-a.

\end{document}